\documentclass{article}
\usepackage{amsmath, amsthm, amsfonts}

\newtheorem{theorem}{Theorem}[section]

\newtheorem{lemma}[theorem]{Lemma}
\newtheorem{Lem}{Lemma}[section]
\newtheorem{definition}{Definition}[section]

      \def\ni{\noindent}
      \def\nn{\nonumber}
      \def\rf#1{\mbox{$(\ref{#1})$}}
      
      \def\l{\ell}

      \def\be{\begin{equation}} 
      \def\ee{\end{equation}} 
      \def\beqn{\begin{eqnarray}} 
      \def\eeqn{\end{eqnarray}} 
      \def\beq{\begin{eqnarray*}} 
      \def\eeq{\end{eqnarray*}}
      \def\proof{{\ni\bf Proof:}\ }
      \def\mb{\mbox} 
      \def\de{\delta}

      \def\ra{\rightarrow} 

      \def\de{\delta}
      \def\lan{\langle}
      \def\ran{\rangle}
      \def\goto{\rightarrow}

\newcommand{\IR}{{\mathbb R}}
\newcommand{\IZ}{{\mathbb Z}}
\def\p{{M_1(E)}}
\def\l{\langle}
\def\r{\rangle}
\newcommand\eproof{\unskip\enskip\null\nobreak\hfill $\sqcup\llap{$\sqcap$}$\par}
\begin{document}

\title{\Large {\bf Reversibility of Interacting Fleming-Viot Processes with Mutation, Selection, and Recombination}
\thanks{Research supported by the Natural Sciences and Engineering Research Council of Canada}}
\author{Shui Feng\\McMaster University \and Byron Schmuland\\ University of Alberta\and Jean
Vaillancourt\\Universit\'e du Qu\'ebec en Outaouais \and Xiaowen
Zhou\\Concordia University}

\maketitle

\begin{abstract}
Reversibility of the Fleming-Viot process with mutation, selection,
and recombination is well understood. In this paper, we study the
reversibility of a system of Fleming-Viot processes that live on a
countable number of colonies interacting with each other through
migrations between the colonies. It is shown that reversibility
fails when both migration and mutation are non-trivial.
\end{abstract}

\section{Introduction}
The Fleming-Viot process is a probability-measure-valued Markov
process describing the evolution of the distribution of allelic
types in a large population. It arises most naturally in population
genetics as the limit in distribution of certain sequences of Markov
chains undergoing mutation, natural selection, recombination, and
random genetic drift.

Reversibility plays an important role in statistical inference in
the neutral theory of population genetics. When reversibility holds,
techniques used for future predictions can then be used to
understand the starting distribution that lead to the present state.
Several models, such as the Wright-Fisher Markov chain and the
finite alleles Wright-Fisher diffusion, are reversible. The
reversibility of the Fleming-Viot process with parent independent
mutation was obtained in \cite{Et90} and \cite{shiga90}. On the
other hand, reversibility is a very restrictive property. The
results in \cite{over-roc97}, \cite{lsy99}, and \cite{handa02}, show
that the Fleming-Viot process is reversible only if the mutation,
natural selection, and recombination have a special form.

The interacting Fleming-Viot process studied in this paper is a
countable collection of Fleming-Viot processes that interact through
geographical migration. It is the diffusion approximation to the
stepping-stone model involving infinitely many alleles. Without
migration, our system will simply be a collection of independent
Fleming-Viot processes. The migration can be viewed as an external
force acting upon the independent system of the Fleming-Viot
processes. Since the internal reversible forces such as mutation and
selection are constantly corrected by the external migration force,
it is natural to expect the loss of reversibility in the interacting
Fleming-Viot process due to competition between local forces and
migration.

The long-time behavior of the interacting Fleming-Viot process is
well-known. In the absence of mutation, selection, and
recombination, a complete characterization of stationary
distributions were obtained in \cite{shiga80} for the two allele
case and in \cite{dgv95} for the general case in terms of migration.
In \cite{shigauchi86} (two allele) and \cite{dg99} (general), the
structures of the stationary distributions were investigated for
models involving mutation, selection and recombination. In this
paper we study the reversibility of the general Fleming-Viot process
and investigate the interrelation between mutation, selection and
recombination, and migration. Under very general hypotheses, we show
that the interacting Fleming-Viot process with mutation, selection,
recombination, and migration is irreversible. Our results cover all
models in \cite{shiga80}, \cite{shigauchi86}, \cite{dgv95}, and
\cite{dg99}.

\section{Model}

      Let $I$ be a countable index set where each element $\xi\in I$ labels a colony.
      The different genetic types of individuals in the population will be
      modelled by a compact metric space $E$.
      Let $M_1(E)$ denote the space of Borel probability measures on $E$,  $M(E)$
      be the space of finite signed Borel measures on $E$, and
      $\Delta$ the set of Dirac measures on $E$.  We let $B(E)$ denote the space of bounded measurable functions
      on $E$, and $C(E)$ the space of continuous functions on $E$. For any $\mu$ in
      $M(E)$ and $g$ in $B(E)$, we use the notation  $\langle \mu, g\rangle=\int_E g(x)\,\mu(dx).$
      Let
      \begin{eqnarray*}
      B(E)^I&:=&\{{\bf f}=(f_{\xi})_{\xi\in I}: f_\xi\in B(E)\}\\
      M(E)^I&:=& \{X=(X_{\xi})_{\xi \in I}: X_\xi\in M(E)\}.
      \end{eqnarray*}
      For $X$ in $M(E)^I$ and $\bf f$ in $B(E)^I$, we write $\langle X, {\bf f}\rangle:=\sum_{\xi \in I}\langle X_{\xi}, f_{\xi} \rangle$
       whenever the sum converges.  The state space for our process will be ${\p}^I\subseteq M(E)^I.$

      For every $\xi, \xi'$ in $I$, let $a(\xi,\xi')$ denote the migration probability from colony $\xi$
      to colony $\xi'$. We assume
      \begin{equation} \label{mig1}
      a(\xi,\xi)=0, \ \sum_{\xi' \in I}a(\xi,\xi')=1.
      \end{equation}

Define the mutation operator $(A,{\cal D}(A))$ to be the generator
of a conservative Feller semigroup $(P_t)$ on $C(E)$. We assume that
the domain ${\cal D}(A)$ of $A$ is dense in $C(E)$.

 The sets $C(E)^I$ and ${\cal D}(A)^I$ denote subsets of $B(E)^I$, where the coordinate functions are in
  $C(E)$ and ${\cal D}(A)$, respectively. Set
  $$B(E)^I_0:= \{{\bf f}\in B(E)^I : f_{\xi}\equiv 0\ \mb{for all $\xi$ outside a finite subset of } I\},$$
  and define $C(E)^I_0$ and ${\cal D}(A)^I_0$ similarly.

For any symmetric bounded measurable function $V$ on $E^2$, we
define the selection operator $S: M_1(E)\to M(E)$ by
    $$
S(\mu)(du):=\left(\int_E V(u,v)\mu(dv)- \int_E \int_E
V(v,w)\mu(dv)\mu(dw)\right)\mu(du).
    $$

When two types $u,v$ undergo recombination; the distribution of the
resulting type is distributed according to the probability kernel
$\eta(u,v;dw)$. The recombination operator $R: M_1(E)\to M(E)$ is
given by
 $$
 R(\mu)(du):= \int_E \int_E \eta(v,w;du)\mu(dv)\mu(dw)-\mu(du).
 $$

 Let ${\tilde {\cal A}}$ be the algebra of functions on ${\p}^I$ given by the collection of linear
 combinations of functions of the form
      \begin{equation}\label{f1}
      F(X)=\prod_{i=1}^m \lan X_{\xi_i}, f_i\ran,
      \end{equation}
 where $m\geq 1$, $f_i \in B(E)$ for $1\leq i\leq m$, and $(\xi_1,\ldots,\xi_m) \in I^m$.
 Similarly, let ${\cal A}$ be the sub-algebra of ${\tilde {\cal A}}$, given by linear
 combinations of functions of the form (\ref{f1}) with $f_i\in {\cal D}(A)$ for $1\leq i\leq m.$
Note that both ${\tilde {\cal A}}$ and ${{\cal A}}$ are measure
determining on ${\p}^I$.

For $F: {\p}^I\to\IR$ we define  partial derivatives as follows,
whenever the limit exists:
  $${\partial F(X)\over\partial X_\xi}(u):=\lim_{\varepsilon\downarrow 0}
  {F(X^\varepsilon(\xi,u))-F(X)\over\varepsilon}\quad \mbox{for }u\in E, \, \xi\in I,$$
  with
  $$ (X^\varepsilon(\xi,u))_{\xi^\prime}:=\begin{cases}X_{\xi'}&\mbox{ if }\xi^\prime\not=\xi,\\
   X_\xi+\varepsilon\delta_u &\mbox{ if }\xi^\prime=\xi.\end{cases}$$

This definition requires us to extend the domain of $F$
infinitesimally from ${\p}^I$ to $M(E)^I$. For $F$ in ${\tilde {\cal
A}}$, this is done via  (\ref{f1}).

For non-negative numbers $s,r,\rho$, the generator ${\cal
L}_{s,r,\rho}$ of the interacting Fleming-Viot process incorporating
migration, mutation, selection, and recombination is defined for
$F\in {\cal A}$ by
$$
{\mathcal L}_{s,r,\rho}F(X):={\mathcal L}_{s,r}F(X)+{\mathcal
L}_{\rho}F(X),
$$
where
\begin{equation}\label{immipart}
{\mathcal L}_{\rho}F(X):=\rho\sum_{\xi,\xi' \in I}a(\xi,\xi')
\left\langle X_{\xi'}-X_{\xi}, \frac{\delta F}{\delta
X_{\xi}(\cdot)}\right\rangle,
\end{equation}
\beqn {\mathcal L}_{s,r}F(X)&:=& \sum_{\xi \in I}\left\langle
X_{\xi}, A\frac{\delta F}{\delta X_{\xi}(\cdot)}\right\rangle
      +s\sum_{\xi \in I}\left\langle S(X_{\xi}), \frac{\delta F}{\delta X_{\xi}(\cdot)}\right\rangle
      +r\sum_{\xi \in I}\left\langle R(X_{\xi}),\frac{\delta F}{\delta X_{\xi}(\cdot)}\right\rangle\nn\\[5pt]
    &&+{1\over2}\sum_{\xi\in I}\int_E\int_E \frac{\delta^2 F}{\delta X_{\xi}(u)\delta
      X_{\xi}(v)}Q_{X_{\xi}}(du,dv),\nn
\eeqn and
\[Q_{\mu}(du,dv):=\mu(du)\delta_u(dv)-\mu(du)\mu(dv).\]
For $X\in {\p}^I$ and ${\bf f}\in{\cal D}(A)_0^I$, define
$$ \langle b_\xi(X), f_\xi\rangle
:=  \langle X_{\xi},Af_{\xi}\rangle +\rho\sum_{\xi'\in I}a(\xi,\xi')
\langle X_{\xi'}-X_{\xi},f_{\xi}\rangle + \langle s S(X_{\xi})+r
R(X_{\xi}),f_{\xi}\rangle, $$ and let $\langle b(X),{\bf f}\rangle
:=\sum_{\xi\in I} \langle b_\xi(X), f_\xi\rangle.$
 The generator ${\cal L}_{s,r,\rho}$ can then be written as
\begin{equation}\label{gradform}
{\mathcal L}_{s,r,\rho}F(X)= \left\langle b(X), \frac{\delta
F}{\delta X}\right\rangle
      + {1\over2}\sum_{\xi\in I}\int_E\int_E  \frac{\delta^2 F}{\delta X_{\xi}(u)\delta X_{\xi}(v)}\,Q_{X_{\xi}}(du,dv),
\end{equation}
where $\frac{\delta F}{\delta X}=\left(\frac{\delta F}{\delta
X_{\xi}}\right)_{\xi\in I}$.

\begin{theorem}
For each $X$ in ${\p}^I$, the martingale problem associated with
generator $({\cal L}_{s,r,\rho}, {\cal A})$ starting at $X$ is
well-posed.
\end{theorem}
\proof The case of $\rho=0$, and the case of $A=0, s=r=0$ can be
found respectively in \cite{EtKu93} and \cite{dgv95}. The case of
$r=0$ was obtained in \cite{handa90}. The general case was studied
in \cite{dg99}, where the index set $I$ is either the finite
dimensional lattice or the hierarchical group, and the type space is
the set of integers.

Even though the index set and state space in our model are more
general, the proofs are similar to that used in \cite{handa90} and
\cite{dg99}. For completeness, we sketch a  proof below.

Following \cite{EtKu93b}, define the following system of
Wright-Fisher type Markov chains. For each colony $\xi$ in $I$,
consider a population of $N$ individuals with types in the space
$E$. The population evolves under the influence of mutation,
selection, recombination, migration, and genetic drift. Future
generations are formed as follows: each individual chooses a pair in
the current generation as parents. The probability that a particular
pair is chosen is weighted by the fitness (described by $V(x,y)$) of
the pair. After the parents are selected, a recombination of the
parent types occurs. The type created through recombination will
change again:  first through migration and then mutation. Existence
for the martingale problem  follows from the tightness of the
empirical processes of approximating systems of Markov chains.

Uniqueness follows from the existence of a dual process. Let
\[
{\cal H}:=\bigcup_{m=1}^{\infty}(C(E^m)\times I^m).
\]
For each solution $X(t)=(X_{\xi}(t))$ to the martingale problem
associated with ${\cal L}_{s,r,\rho}$, the law of $X(t)$ is
determined by
\[
F((f,\pi),X(t))=E^{X(0)}\left( \int_E \cdots \int_E f(u_1,\ldots,
u_m)X_{\xi_1}(du_1)\cdots X_{\xi_m}(du_m)\right)
\]
for all $(f, \pi)$ in $C(E^m)\times I^m, m\geq 1.$

For $F(X)=\prod_{i=1}^m \lan X_{\xi_i}, f_i\ran$ in ${\cal A}$,
direct calculations give \beqn \label{duality1} {\mathcal
L}_{s,r,\rho}F(X)&=&\sum_{i=1}^m \left\{\langle X_{\xi_i}, Af_i
\rangle
     +s \langle S(X_{\xi_i})+r R(X_{\xi_i}), f_i\rangle
     + \rho\sum_{\xi' \in I}a(\xi_i,\xi') \langle X_{\xi'}-X_{\xi_i}, f_i\rangle\right\}\prod_{j\neq i}\langle X_{\xi_j}, f_j\rangle\nn\\[5pt]
    &&+ \sum_{1\leq i<k\leq m, \xi_i=\xi_k}(\langle X_{\xi_i},f_i f_k\rangle-\langle X_{\xi_i},f_i\rangle \langle X_{\xi_k}, f_k\rangle)\prod_{j\neq i,k}\langle X_{\xi_j}, f_j\rangle.
\eeqn

Define for $\pi=(\xi_1,\dots,\xi_m)$ in $I^m$, $m \geq 1$ and
$f(u_1,\ldots,u_m)=\prod_{i=1}^m f_i(u_i)$ \beq
&& X_{\pi}(du_1,\ldots,du_m):=\prod_{i=1}^m X_{\xi_i}(du_i),\\
&& \tilde{\pi}^i:=(\xi_1,\ldots,\xi_m,\xi_i), i=1,\ldots,m,\\
&& \tilde{\pi}^{ii}:=(\xi_1,\ldots,\xi_m,\xi_i,\xi_i), i=1,\ldots,m,\\
&& \hat{\pi}^j := (\xi_1,\ldots,\xi_{j-1},\xi_{j+1},\ldots,\xi_m), j=2,\ldots,m,\\
&& \pi^{i,\xi}:=(\xi_1,\ldots,\xi_{i-1},\xi,\xi_{i+1},\ldots,\xi_m),
\eeq and \beq
&& A^m f(u_1,\ldots,u_m) := \sum_{i=1}^m Af_i(u_i) \prod_{j\neq i} f_j(u_j),\\
&& H_{im}f(u_1,\ldots,u_m):= (V(u_i,u_m)-V(u_m,u_{m+1}))f(u_1,\ldots,u_m),\\
&& K_{i}f(u_1,\ldots,u_m):=\int_E f(u_1,\dots,\nu,
u_{i+1},\ldots,u_m)\eta(u_i,u_{m+1};d\nu)-f(u_1,\ldots,u_m). \eeq
Then \rf{duality1} can be written as
\begin{eqnarray*}
{\mathcal L}_{s,r,\rho}F(X)&=& \langle X_{\pi}, A^m f\rangle \\[5pt]
                           &&+\sum_{i=1}^m \left\{s\langle X_{\tilde{\pi}^{ii}}, H_{im}f\rangle
                            +r\langle X_{\tilde{\pi}^i}, K_i f\rangle
                            + \rho \sum_{\xi' \in I}a(\xi_i,\xi') \langle X_{\pi^{i,\xi'}}-X_{\pi}, f\rangle\right\}\\
                           &&+\sum_{1\leq i<k\leq m}(\langle X_{\hat{\pi}^k},\Phi_{ik}f\rangle-\langle X_{\pi},f\rangle),
\end{eqnarray*}
where $\Phi_{ik}f$ is the function in $B(E^{m-1})$ that is obtained
from $f$ by replacing $u_k$ with $u_i$ and relabeling the variables.

The dual process $(f_t,\pi_t)$ is an ${\cal H}$-valued process,
starting from $(f_0,\pi_0)=(f,\pi)$,
    that involves the following transitions:
\begin{itemize}
\item coordinates of $\pi_t$ are independent continuous time Markov chains on $I$ with transition rate $(\rho a(\xi,\xi'))_{\xi,\xi' \in I}$;
\item any two coordinates of $\pi_t$ that are the same will coalesce into one element at the same site with rate one;
\item at rate $s$ a coordinate of $\pi_t$ will create two copies of itself so that the size of $\pi_t$ is increased by two;
\item at rate $r$ a coordinate of $\pi_t$ will create a copy of itself so that the size of $\pi_t$ is increased by one;
\item $f_0$ is in $C(E^{|\pi_0|})$; between transitions of $\pi_t$, $f_t$ follows a deterministic path
determined by the semigroup associated with $|\pi_t|$ independent
copies of $A$-motion;
\item At the time of coalescence, the corresponding variables in $f_t$ are set equal, which results in a jump from space
$C(E^{|\pi_{t-}|})$ to space $C(E^{|\pi_{t-}|-1})$;
\item If two new coordinates are created when the current number of variables is $m$, then we have
\[
f(u_1,\ldots,u_m) \ra
s(V(u_i,u_{m+1})-V(u_{m+1},u_{m+2}))f(u_1,\ldots,u_m);
\]
\item If one new coordinate is created when the current number of variables is $m$, then we have
\[
f(u_1,\ldots,u_m) \ra \int_E f(u_1,\ldots,u_{i-1},\nu
,u_{i+1}\ldots,u_m)\eta(u_i,u_{m+1};d\nu).
\]
\end{itemize}

The uniqueness now follows from the following duality relation
\[
E_{X(0)}[\langle X_{\pi}(t), f\rangle]=E_{(f,\pi)}[\langle
X_{\pi_t}(0), f_t\rangle\ e^{s\int_0^t |\pi_u|du}\rangle].
\]
\eproof

\section{ Quasi-invariance and the cocycle identity} \setcounter{equation}{0}

In this section we prove the main result of the paper relating the
reversibility of probability measures on $\p$ with their
quasi-invariance. These results generalize those proved by Handa for
the single site Fleming-Viot process. In the sections that follow,
we will show that reversibility is a very restrictive condition that
only applies to very special cases of the Fleming-Viot model.

\begin{definition}\label{reversible}
A probability measure $\Pi$ on $\p$ is reversible for the
Fleming-Viot operator $({\cal L}_{s,r,\rho}, {\cal A})$ if  for
$\Phi,\Psi\in{\cal A}$,
\[ \int {\cal L}_{s,r,\rho}\Phi(X) \Psi(X)\,\Pi(dX)= \int {\cal L}_{s,r,\rho}\Psi(X) \Phi(X)\,\Pi(dX). \]
\end{definition}

For each ${\bf f}$ in $C(E)^I$, define a map $S_{{\bf f}}: {\p}^I
\to {\p}^I$ by $ S_{{\bf f}}(X)= (X^{f_{\xi}}_{\xi})_{\xi \in I}$,
where
\[ X_{\xi}^{f_{\xi}}(dv):={e^{f_{\xi}(v)}X_{\xi}(dv)\over \langle X_\xi, e^{f_\xi} \rangle}. \]
It follows from the definition that $S_{\bf f}(S_{\bf g})= S_{{\bf
f}+{\bf g}}$ for any ${\bf f}, {\bf g}$ in $C(E)^I$. For any ${\bf
f}$ in $C(E)^I$ and probability measure $\Pi$ on $\p^I$, set
$\Pi^{{\bf f}}(\cdot):=\Pi(S_{{\bf f}}(\cdot))$.

The probability $\Pi$ is called {\it quasi-invariant} for ${\cal
D}(A)^I_0$ if for any ${\bf f}\in {\cal D}(A)^I_0$,  the measures
$\Pi^{{\bf f}}$ and $\Pi$ are mutually absolutely continuous with
\[
\frac{d \Pi^{{\bf f}}}{d\Pi}(X)=\exp\{\Lambda({\bf f}, X)\},
\]
where $\Lambda: {\cal D}(A)^I_0\times {\p}^I \mapsto \IR$ is called
the {\it cocycle} associated with $\Pi$.

A direct result of the quasi-invariance is the following  {\it
cocycle identity}: for any ${\bf f}, {\bf g} \in {\cal D}(A)^I_0$,
for $\Pi$ almost all $X$,
\begin{equation} \label{cocyle}
 \Lambda({\bf f}+{\bf g},X)= \Lambda({\bf f},S_{\bf g}(X))+\Lambda({\bf g},X).
\end{equation}

The {\it carr\'e du champ} associated with the operator ${\mathcal
L}_{s,r,\rho}$ is defined by
\begin{equation}\label{def-cdc}
\Gamma(\Phi,\Psi)={1\over2} ({\mathcal L}_{s,r,\rho}(\Phi\Psi) -\Phi
{\mathcal L}_{s,r,\rho}(\Psi)-{\mathcal
L}_{s,r,\rho}(\Phi)\Psi),\quad \Phi,\Psi\in{\cal A}.
\end{equation}


For any two functions $f, g$ in $B(E)$, set $(f\otimes
g)(u,v):=f(u)g(v)$. By an argument similar to that used in the
proof of Lemma 3.1 in \cite {handa02}, we obtain the following
result.

\begin{Lem}\label{le1} For $\Phi,\Psi\in{{\cal A}}$ and $X\in {\p}^I$,
\begin{equation}\label{le11}
\Gamma(\Phi,\Psi)(X)={1\over 2}\sum_{\xi\in I}\left\l Q_{X_\xi},
{{\de \Phi(X)\over\de X_\xi}}\otimes{{\de \Psi(X)\over\de
X_\xi}}\right\r,
\end{equation}
and for $\Phi, \Psi_1,\Psi_2\in{{\cal A}}$,
\begin{equation}\label{le12}
\Gamma(\Phi\Psi_1,\Psi_2)+\Gamma(\Phi\Psi_2,\Psi_1)-\Gamma(\Phi,\Psi_1\Psi_2)=2\Phi\Gamma(\Psi_1,\Psi_2).
\end{equation}
\end{Lem}

\begin{Lem}\label{le1b}
The probability measure $\Pi$ is reversible with respect to ${\cal
L}_{s,r,\rho}$ if and only if
\begin{equation}\label{le1b1}
-\frac{1}{2}\int \left\lan  Q_{X_\xi},\frac{\de \Phi(X)}{\de X_\xi}
\otimes f_\xi \right\ran\Pi(dX)= \int \Phi(X)\lan
b_{\xi}(X),f_\xi\ran\Pi(dX)
\end{equation}
for any $\Phi\in{{\cal A}}$, $\xi\in I$, and $f_\xi\in {\cal D}(A)$.
\end{Lem}

\proof Assume that $\Pi$ is reversible with respect to ${\cal
L}_{s,r,\rho}$. For a fixed $\xi$ in $I$, let $\Psi(X)=\langle
X_{\xi}, f_{\xi}\rangle$. It follows from \rf{gradform} that  ${\cal
L}_{s,r,\rho}\Psi(X)= \langle b_{\xi}(X),f_{\xi}\rangle$. This,
combined with Lemma \ref{le1} and reversibility, implies
(\ref{le1b1}).

Next we assume that (\ref{le1b1}) holds.  First we show, by
induction on $n$, that for any $n\geq 1$
\begin{equation}\label{le1b2}
\int \Phi(X){\mathcal
L}_{s,r,\rho}\Psi^{(n)}(X)\,\Pi(dX)=-\int\Gamma(\Phi,\Psi^{(n)})(X)\,\Pi(dX),
\end{equation}
for any $\Phi\in{{\cal A}},  f_{i}\in {\cal D}(A), \xi_i\in I,
i=1,\ldots,n$, and
\[\Psi^{(n)}(X):=\prod_{i=1}^n\Psi_i(X):=\prod_{i=1}^n\lan X_{\xi_i}, f_i \ran.\]

The case of $n=1$ follows from (\ref{le11}) and (\ref{le1b1}).
Assume that (\ref{le1b2}) holds for $n\leq k$. It follows from
(\ref{def-cdc}) and (\ref{le12}) that
\begin{equation}
\begin{split}
\Phi {\mathcal L}_{s,r,\rho}(\Psi^{(k+1)})&=\Phi[2\Gamma(\Psi^{(k)},
\Psi_{k+1})
+ \Psi_{k+1}{\mathcal L}_{s,r,\rho}(\Psi^{(k)})+\Psi^{(k)} {\mathcal L}_{s,r,\rho}\Psi_{k+1}]\\
&=\Gamma(\Phi\Psi^{(k)},\Psi_{k+1})+\Phi\Psi^{(k)}
{\mathcal L}_{s,r,\rho}\Psi_{k+1}\\
&\quad + \Gamma(\Phi\Psi_{k+1},\Psi^{(k)})+ \Phi\Psi_{k+1}{\mathcal L}_{s,r,\rho}(\Psi^{(k)})\\
&\quad-\Gamma(\Phi,\Psi^{(k+1)})
\end{split}
\end{equation}
which implies that
\[\int \Phi(X) {\mathcal L}_{s,r,\rho}(\Psi^{(k+1)})(X) \Pi(dX)=-\int \Gamma(\Phi,\Psi^{(k+1)})(X)\Pi(dX).\]

It follows from (\ref{le1b2}) that for any $ \Phi, \Psi$ in ${{\cal
A}}$
\[\int \Phi(X){\mathcal L}_{s,r,\rho}\Psi(X)\Pi(dX)=-\int \Gamma(\Psi,\Phi)(X)\Pi(dX);\]
and by symmetry
\[\int \Psi(X){\mathcal L}_{s,r,\rho}\Phi(X)\Pi(dX)= \int \Phi(X){\mathcal L}_{s,r,\rho}\Psi(X)\Pi(dX). \]
Therefore, $\Pi$ is reversible with respect to ${\cal
L}_{s,r,\rho}$. \eproof

\begin{Lem}\label{le2} Suppose ${\bf f}\in C(E)^I$ and put $X_t:=S_{-t{\bf f}}X$ for $X\in {\p}^I$ and $t\in\IR$.
For every $\Phi\in{\tilde {\cal A}}$ we have
\begin{equation}\label{need}
{d\over dt}\Phi(X_t)=-\sum_{\xi\in I}\left\l Q_{X_\xi},
f_\xi\otimes{{\de \Phi(X_t)\over\de X_\xi}} \right\r.
\end{equation}
\end{Lem}

\proof  Since both sides of the equation are linear, it suffices to
prove the result for functions of the form $\Phi(X)= \prod_{i=1}^m
\langle X_{\xi_i}, g_i\rangle$, where $m$ a positive integer,
$(\xi_i)_{1\leq i\leq m}$ in $I$, and $g_i\in B(E)$. But
 both sides of the equation are also derivations in $\Phi$, so it suffices
to take $m=1$. But in this case, (\ref{need}) follows from an easy
calculation or Lemma 3.3 of \cite{handa02}. \eproof

For ${\bf f}\in {\cal D}(A)^I_0$ and $X\in {\p}^I$, we let
\begin{equation}\label{cocycle-def}
\begin{split}
\Lambda({\bf f},X)&:=2\int^1_0 \l b(S_{s{\bf f}}X),{\bf
f}\rangle\,ds.
\end{split}
\end{equation}

\begin{Lem}\label{le3} Suppose ${\bf f}\in {\cal D}(A)^I_0$, and put $X_t=S_{-t{\bf f}}X$ for $X\in {\p}^I$ and $t\in\IR$.
Then we can write
\[\Lambda(t{\bf f},X_t)=2\int_0^t\l b(X_s),{\bf f}\r\,ds.\]
\end{Lem}
\proof
\begin{eqnarray*}
\Lambda(t{\bf f},X_t)&=& 2\int_0^1 \l b(S_{st{\bf f}}X_t),t{\bf f}\r\,ds = 2t\int_0^1 \l b(S_{-(1-s)t{\bf f}}X),{\bf f}\r\,ds  \\
&=&  2t\int_0^1 \l b(S_{-st{\bf f}}X),{\bf f}\r\,ds = 2\int_0^t \l
b(S_{-s{\bf f}}X),{\bf f}\r\,ds.
\end{eqnarray*}   \eproof

The following lemma proves formula (\ref{need}) for certain
functions $F\notin \tilde{{\cal A}}$.


\begin{Lem}\label{le4} Suppose ${\bf f}\in C(E)^I_0$, and put $X_t=S_{-t{\bf f}}X$ for $X\in {\p}^I$ and $t\in\IR$.
For $h\in C(E)$ and the sequence $c(\xi)$ satisfying $\sum_{\xi\in
I}|c(\xi)|<\infty$, define $F:{\p}^I\to\IR$ by \[F(X):=\left\l
\sum_{\xi\in I}c(\xi) X_\xi, h \right\r.\] Then
$$
{d\over dt}F(X_t)=-\sum_{\xi\in I}\left\l Q_{X_\xi}, \,
f_\xi\otimes{{\de F(X_t)\over\de X_\xi}} \right\r.
$$
\end{Lem}
\proof Let $I_0$ be a finite subset of $I$ such that $f_\xi=0$ for
$\xi\not\in I_0$. Define
\[
F_0(X):=\left\l \sum_{\xi\in I_0}c(\xi) X_\xi, h\right\r.
\]
Clearly $F_0\in {\tilde {\cal A}}$. Also, $(X_t)_\xi=X_\xi$ for
$\xi\notin I_0$ so those terms have a zero time derivative.
Therefore,
$${d\over dt}F(X_t)={d\over dt}F_0(X_t).$$
  It follows from direct calculation that
  $$ {{\de F_0(X_t)\over\de X_\xi}} =\begin{cases}{\frac{\de F(X_t)}{\de X_\xi}}&\text{ \, if \,} \xi\in I_0,\cr
  \noalign{\vskip8pt}    0 & \text{\, if \,} \xi\notin I_0.\end{cases}$$
  Since $f_\xi\equiv 0$ for $\xi\notin I_0$, this gives
  \begin{equation*}
\begin{split}
\sum_{\xi\in I}\left\l Q_{X_\xi}, \, f_\xi\otimes{{\de
F(X_t)\over\de X_\xi}} \right\r
  &= \sum_{\xi\in I_0}\left\l Q_{X_\xi}, \, f_\xi\otimes\frac{\de F(X_t)}{\de X_\xi} \right\r\\
  &= \sum_{\xi\in I_0}\left\l Q_{X_\xi}, \, f_\xi\otimes\frac{\de F_0(X_t)}{\de X_\xi} \right\r\\
  &= \sum_{\xi\in I}\left\l Q_{X_\xi}, \, f_\xi\otimes\frac{\de F_0(X_t)}{\de X_\xi}  \right\r,
\end{split}
\end{equation*}
which, combined with Lemma \ref{le2}, implies the result. \eproof

\begin{theorem}\label{t2}
If the probability measure $\Pi$ in $M_1({\p}^I)$ is reversible for
${\mathcal L}_{s,r,\rho}$, then $\Pi$ is quasi-invariant for ${\cal
D}(A)^I_0$ with cocycle $\Lambda({\bf f}, X)$ given by
\rf{cocycle-def}.
\end{theorem}

\proof Assume that $\Pi\in M_1({\p}^I)$ is reversible with respect
to ${\cal L}_{s,r,\rho}$, and fix ${\bf f}\in {\cal D}(A)^I_0$. We
must show that
$$\int F(X)\,(S_{\bf f}\Pi)(dX)=\int F(S_{-{\bf f}}X)\,\Pi(dX)=\int F(X)e^{\Lambda({\bf f},X)}\,\Pi(dX),$$
for sufficiently many functions $F:{\p}^I\to\IR.$ Since
$\exp(-\Lambda({\bf f},X))$ is strictly positive and ${\cal A}$ is
measure determining, it suffices to prove that for any $\Phi\in
{{\cal A}}$
$$\int \Phi(S_{-{\bf f}}X)e^{-\Lambda({\bf f},S_{-{\bf f}}X)}\,\Pi(dX)= \int \Phi(X)\,\Pi(dX).$$
In what follows we shall show that
$$Z(t):= \int \Phi(S_{-t{\bf f}}X)e^{-\Lambda(t{\bf f},S_{-t{\bf f}}X)}\,\Pi(dX)$$
is a constant function of $t\in\IR$. Setting
$$\tilde\Phi_t(X):=\Phi(X_t)e^{-\Lambda(t{\bf f},X_t)}= \Phi(S_{-t{\bf f}}X)e^{-\Lambda(t{\bf f},S_{-t{\bf f}}X)},$$
and noting that $\Lambda(t{\bf f},X_t)=2\int_0^t\l b(X_s),{\bf
f}\r\,ds$, we have
\begin{equation}\label{t2a}
{\de \tilde\Phi_t(X)\over\de X_\xi}(u) = {\de \Phi(X_t)\over\de
X_\xi}(u) e^{-\Lambda(t{\bf f},X_t)}- 2\tilde\Phi_t(X)\int_0^t {\de
\l b(X_s),{\bf f}\r\over\de X_\xi}(u)\,ds.
\end{equation}
It follows that
\begin{equation*}
\begin{split}
&\sum_{\xi\in I} \left\l Q_{X_\xi}, \, f_\xi\otimes \frac{\de \tilde\Phi_t(X)} {\de X_\xi} \right\r \\
&= \sum_{\xi\in I}\left\l Q_{X_\xi}, \, f_\xi\otimes \frac{\de
\Phi(X_t) }{\de X_\xi} \right\r e^{-\Lambda(t{\bf f},X_t)}
  - 2\tilde\Phi_t(X) \int^t_0 \sum_{\xi\in I}\left\l Q_{X_\xi}, \,   f_\xi\otimes \frac{\de\l b(X_s),{\bf f}\r} {\de X_\xi} \right\r ds \\
  &= -{d\over dt} \Phi(X_t)\, e^{-\Lambda(t{\bf f},X_t)} + 2\tilde\Phi_t(X) \int^t_0 {d\over ds} \l b(X_s),{\bf f}\r \,ds \\
  &= -{d\over dt} \Phi(X_t)\, e^{-\Lambda(t{\bf f},X_t)} + 2\tilde\Phi_t(X)(\l b(X_t),{\bf f}\r -\l b(X),{\bf
  f}\r),
\end{split}
  \end{equation*}
where Lemmas \ref{le2} and \ref{le4} are used for obtaining the
second equality. Therefore,
  \begin{equation*}
\begin{split}
&Z^\prime(t)\\
&=\int \left({d\over dt}\Phi(X_t)\, e^{-\Lambda(t{\bf f},X_t)}+\Phi(X_t){d\over dt} e^{-\Lambda(t{\bf f},X_t)} \right)\Pi(dX) \\
&=\int \left( -\sum_{\xi\in I} \left\l Q_{X_\xi}, \, f_\xi\otimes
{\de \tilde\Phi_t(X)\over\de X_\xi} \right\r
          +2\tilde\Phi_t(X) (\l b(X_t),{\bf f}\r -\l b(X),{\bf f}\r) -2\tilde\Phi_t(X) \l b(X_t),{\bf f}\r \right)\Pi(dX) \\
&=-\int  \sum_{\xi\in I} \left\l Q_{X_\xi}, \, f_\xi\otimes {\de
\tilde\Phi_t(X)\over\de X_\xi} \right\r \,\Pi(dX)
          -2\int \tilde\Phi_t(X) \l b(X), {\bf f}\r \,\Pi(dX).
\end{split}
  \end{equation*}

 By reversibility and Lemma \ref{le1b},
\begin{equation}\label{t2b}
\int \sum_{\xi\in I} \left\l Q_{X_\xi}, \, f_\xi \otimes
{\Phi(X)\over\de X_\xi} \right\r
  \,\Pi(dX)  + 2\int \l b(X), {\bf f}\r\Phi(X)  \,\Pi(dX) =0,
\end{equation}
  for $\Phi\in{{\cal A}}$. In the Appendix, we introduce a space of functions ${\mathcal H}$ that contains ${\mathcal A}$,
  and show that $ \tilde\Phi_t(X)\in {\cal H}$ and \rf{t2b} holds for all $\Phi$ in ${\mathcal H}$.
  These implie that $Z^\prime(t)=0$. Therefore, $Z(1)= Z(0)$ and the theorem follows from
  \[  \int \Phi(S_{-{\bf f}}X)e^{-\Lambda({\bf f},S_{-{\bf  f}}X)}\,\Pi(dX)=Z(1)= Z(0)= \int \Phi(X)\Pi(dX). \]  \eproof

\begin{theorem}\label{t3}
If the probability measure $\Pi$ in $M_1({\p}^I)$ is quasi-invariant
with cocycle given by (\ref{cocycle-def}), then $\Pi$ is reversible
with respect to ${\mathcal L}_{s,r,\rho}$.
\end{theorem}

\proof Suppose that $\Pi$ is quasi-invariant with cocycle given by
(\ref{cocycle-def}). Then for any $\xi \in I$ and ${\bf f}$ in
$C(E)^I$ such that $f_{\xi} \in {\cal D}(A)$ and $f_{\xi'}=0$ for
$\xi'\neq \xi$, the function
$$Z(t)= \int \Phi(S_{-t{\bf f}}X)e^{-\Lambda(t{\bf f},S_{-t{\bf f}}X)}\,\Pi(dX)$$
is constant in $t\in\IR$. Noting that
\begin{equation*}
\begin{split}
0=Z^\prime(0) &=-\int  \sum_{\xi\in I} \left\l f_\xi\otimes {\de
\Phi(X)\over\de X_\xi},Q_{X_\xi} \right\r \,\Pi(dX)
          -2\int \Phi(X) \l b(X),{\bf f}\r \,\Pi(dX), \\
\end{split}
\end{equation*}

and $f_{\xi}$ is arbitrary in ${\cal D}(A)$, the theorem follows
from Lemma \ref{le1b}. \eproof

\section{Consequences of the cocycle identity} \setcounter{equation}{0}

    It follows from the cocycle identity \rf{cocyle} that for any $X$ in ${\p}^I$ and any ${\bf f}, {\bf g} \in {\cal D}(A)^I_0$,
\begin{equation}\label{cocyid1}
\Lambda({\bf f},S_{\bf g}(X))-\Lambda({\bf f},X)=\Lambda({\bf
g},S_{\bf f}(X))-\Lambda({\bf g},X).
    \end{equation}

      For any two distinct $\xi_1,\xi_2$ in $I$, and $f,g$ in ${\cal D}(A)$, let
      ${\bf f}=(f_{\xi})$ and ${\bf g}=(g_{\xi})$ be such that $f_{\xi_1}=f, f_{\xi}=0$ for $\xi\neq \xi_1$, and
      $g_{\xi_2}=g,g_{\xi}=0$ for $\xi\neq \xi_2$. By direct calculation,
    \beq
      \Lambda({\bf f},X)&=&2\int_0^1\{\langle S_{u{\bf f}}(X)_{\xi_1},Af\rangle +s \langle
S(S_{u{\bf f}}(X)_{\xi_1}), f\rangle\\
&& + r\langle R(S_{u{\bf f}}(X)_{\xi_1}), f \rangle +
\rho\sum_{\xi'}a(\xi_1,\xi')
    \langle S_{u{\bf f}}(X)_{\xi'}-S_{u{\bf f}}(X)_{\xi_1},f\rangle\}du\\
&=&2\int_0^1\{\langle X^{uf}_{\xi_1},Af\rangle +s \langle
S(X^{uf}_{\xi_1}), f\rangle\\
&& + r\langle R(X^{uf}_{\xi_1}), f \rangle + \rho\sum_{\xi'\neq
\xi_1}a(\xi_1,\xi')
    \langle X_{\xi'}-X^{uf}_{\xi_1},f\rangle\}du,
\eeq and \beq
      \Lambda({\bf f},S_{\bf g}(X))&=&2\int_0^1\{\langle S_{u{\bf f}+{\bf g}}(X)_{\xi_1},Af\rangle +s \langle S(S_{u{\bf f}+{\bf g}}(X)_{\xi_1}), f \rangle\\
&& + r \langle R(S_{u{\bf f}+{\bf g}}(X)_{\xi_1}), f \rangle +
\rho\sum_{\xi'}a(\xi_1,\xi')
    \langle S_{u{\bf f}+{\bf g}}(X)_{\xi'}-S_{u{\bf f}+{\bf g}}(X)_{\xi_1},f\rangle\}du\\
    &=&2\int_0^1\{\langle X^{uf}_{\xi_1},Af\rangle +s \langle S(X^{uf}_{\xi_1}), f \rangle
+ r \langle R(X^{uf}_{\xi_1}), f \rangle\\
&&+ \rho\sum_{\xi'\neq\xi_1,\xi_2}a(\xi_1,\xi')
    \langle X_{\xi'}-X^{uf}_{\xi_1},f\rangle
    +\rho a(\xi_1,\xi_2)\langle
    X^{g}_{\xi_2}-X^{uf}_{\xi_1},f\rangle\}du,
    \eeq
    which leads to

    \begin{equation}\label{cocyid3}
    \Lambda({\bf f},S_{\bf g}(X))-\Lambda({\bf f},X)= 2\rho a(\xi_1,\xi_2)\langle X^{g}_{\xi_2}-X_{\xi_2},f\rangle.
    \end{equation}

    Together, \rf{cocyid1} and \rf{cocyid3} implies that for $\rho >0$
    \begin{equation}\label{cocyid4}
a(\xi_1,\xi_2)\langle
X^{g}_{\xi_2}-X_{\xi_2},f\rangle=a(\xi_2,\xi_1)\langle
X^{f}_{\xi_1}-X_{\xi_1},g\rangle.
    \end{equation}

Let $$\hat{I}:=\{\xi \in I: \mb{there exists }\ \eta\in I,\
\mb{such that}\ a(\eta,\xi)>0\}.$$ It follows from \rf{mig1} that
the set $\hat{I}$ is not empty.

    \begin{lemma}\label{le5}
Suppose that $\Pi$ is a reversible probability measure with respect
to ${\mathcal L}_{s,r,\rho}$ with $\rho>0$.
    Then for any $\xi\in \hat{I}$,  $X_{\xi}$ is a Dirac measure with $\Pi$ probability one, i.e,
      \begin{equation}\label{dirac}\Pi\{X_\xi\in\Delta\}=1.\end{equation}
    \end{lemma}
    \proof Let $C$ be a countable dense subset of $E$.
    By definition, for each $\xi$ in $\hat{I}$, there exists $\xi'$ in $I$ such
    that $a(\xi',\xi)>0$.  Assume that with positive $\Pi$ probability, $X_{\xi}$ is not a Dirac measure.
For any two distinct elements $c_1, c_2$ in $C$, and any positive
rational numbers $r_1, r_2$ satisfying $r_1+r_2<d(c_1,c_2)$, let
\[D(c_1,c_2;r_1,r_2):=\{X\in {\p}^I: X_{\xi}(B(c_1,r_1))>0, X_{\xi}(B(c_2,r_2))>0\},\]
where $B(c_i,r_i)$ denotes the open ball in $E$ with center $c_i$
and radius $r_i$. Clearly,
\[\bigcup_{c_1,c_2;r_1,r_2}D(c_1,c_2,r_1,r_2)=\{X\in{\p}^I: X_\xi\neq \delta_u, \forall  u\in E\}.\]
Therefore, we can find rational numbers $c_1,c_2,r_1,r_2$ such that
$\Pi\{D(c_1,c_2,r_1,r_2)\}>0.$ Choose a nonnegative continuous
function $f$ such that $f(x)=0$ for $x\in B(c_1,r_1)$ and $f(x)>0$
for $x\in B(c_2,r_2)$. For any $X_{\xi}\in D(c_1,c_2;r_1,r_2)$,
observe that $\lan X_{\xi},e^f \ran>1 $.  When the signed measure
$X_{\xi}-X^f_{\xi}$ is restricted to set $B(c_1,r_1)$, we have
\[X_{\xi}-X^f_{\xi}=(1-\lan X_{\xi},e^f \ran^{-1})X_{\xi},\]
which is a measure on $B(c_1,r_1)$ with strictly positive total
mass. Let $g$ be any continuous function such that $g(x)>0$ for
$x\in B(c_1,r_1)$ and $g(x)=0$ for $x\not\in B(c_1,r_1)$.

For any $h\in C(E)$ and any positive integer $k$, define
\[h^{(k)}:=k\int_0^{\frac{1}{k}} P_s h ds.\]
Then $\|h^{(k)}-h\|_\infty\goto 0$, $h^{(k)}\in {\cal D}(A)$ and
$Ah^{(k)}=k(P_{{1}/{k}}h-h)\in C(E)$.

By dominated convergence theorem, we have
\begin{equation*}
\lim_{k\goto\infty}\lan X_{\xi}-X^{f^{(k)}}_{\xi}, ng^{(k)} \ran
=\lan X_{\xi}-X^f_{\xi}, ng \ran = \lan X_{\xi}-X^f_{\xi},
ng1_{B(c_1,r_1)}\ran
\end{equation*}
and
\begin{equation*}
\lim_{k\goto\infty}\lan X_{\xi'}-X^{ng^{(k)}}_{\xi'}, f^{(k)}\ran
=\lan X_{\xi'}-X^{ng}_{\xi'}, f\ran \leq \| f\|_\infty
\end{equation*}
for all $n$.  Choosing $\xi'=\xi_2, \xi=\xi_1$, $g=ng^{(k)}$ in
(\ref{cocyid4}), and taking the limit in the order of $k\ra\infty$
and $n \ra \infty$, gives a contradiction.  \eproof

\noindent {\bf Remark.} It follows from the above theorem that for
each $\xi $ in $\hat{I}$, there is a random variable $x_{\xi}$
taking values in $E$ such that $X_{\xi}=\delta_{x_{\xi}}$ almost
surely under $\Pi$.

\begin{lemma}\label{le6}
Suppose that $\Pi$ is a reversible measure with respect to
${\mathcal L}_{s,r,\rho}$. For each $\xi$ in $I$, let
$I_{\xi}=\{\xi' \in I: a(\xi,\xi')>0\}$. Then for $\xi\in\hat{I}$,
we have $\Pi\{x_\xi= x_{\xi'}\}=1$, for all $\xi'\in I_{\xi}$.
\end{lemma}
\proof By Lemma \ref{le5}, for $\Pi$ almost all $X$, we have
$X_{\xi'}=\delta_{x_{\xi'}}$ for any $\xi'\in I_\xi$. For
$f,g\in{\cal D}(A)$, set $\Phi(X)=\langle X_\xi,f\rangle$ and
$\Psi(X)=\langle X_\xi,g\rangle$. The reversibility, combined with
Lemma \ref{le1}, implies
\begin{equation}\label{equ1}
-\int \Psi(X) {\cal L}_{s,r,\rho}\Phi(X) \Pi(dX)=
\frac{1}{2}\int\lan f\otimes g , Q_{X_\xi}\ran \Pi(dX)=0,
\end{equation}
since $Q_{X_\xi}$ is the zero measure when $X_\xi$ is a delta mass.

Then for any $f, g\in\cal{D}(A)$ equation (\ref{equ1}) gives
\begin{equation}\label{equ2}
\int g(x_\xi)\left[\sum_{\xi'\neq \xi}\rho
a(\xi,\xi')(f(x_{\xi'})-f(x_{\xi}))+Af(x_\xi)+\tilde{R}f(x_\xi)\right]\,\Pi(dX)=0,
\end{equation}
where
\[\tilde{R}f(x):= r\left[\int f(u)\eta(x,x;du)-f(x)\right], \quad x\in E.\]

For any $c\in E$ and $0<r<r'$, choose a sequence of continuous
functions $(f_m)$ on $E$ such that  $0\leq f_m\leq 1$ and $f_m(x)=1$
for $x\in \bar{B}(c,r')$ and $f_m$ converges, pointwisely, to
$1_{\bar{B}(c,r')}$, where $\bar{B}(c,r')$ denotes the closed ball
with center $c$ and radius $r'$; also choose a sequence of
continuous functions $(g_n)$ on $E$ such that $0\leq g_n\leq 1$,
$g_n(x)=1 $ for $x\in \bar{B}(c,r)$, $g_n$ has its support in
$\bar{B}(c,r')$, and $g_n$ converges pointwise to
$1_{\bar{B}(c,r)}$.

By the maximal principle for $A$, we have $Af_m^{(k)}(x)\leq 0$ for
$x\in\bar{B}(c,r')$, so that for $m, n, k, k'$,
\begin{equation}\label{equ3}
\int g_n^{(k')}(x_\xi)Af_m^{(k)}(x_\xi)\,\Pi(dX)\leq 0.
\end{equation}

Since $g^{(k')}_n$ converges pointwise to $g_n $ as $k'\goto \infty$
and $f^{(k)}_m$ converges pointwise to  $f_m $ as $k\goto \infty$,
taking limits in the order of $k'\goto\infty$, $k\goto\infty$,
$m\goto\infty$, and $n\goto\infty$, we first have
\begin{equation}\label{equ3a}
\int g^{(k')}_n(x_\xi)\tilde{R}f^{(k)}_m(x_\xi)\,\Pi(dX)\goto
r\int1_{\bar{B}(c,r)}(x_\xi)
\left(\eta(x_\xi,x_\xi;\bar{B}(c,r'))-1_{\bar{B}(c,r')}(x_{\xi})\right)\,\Pi(dX)\leq
0,
\end{equation}
then combining (\ref{equ2}), (\ref{equ3}) and (\ref{equ3a}) we
further have
\begin{equation}\label{equ5}
\int 1_{\bar{B}(c,r)}(x_\xi)\sum_{\xi'\neq
\xi}a(\xi,\xi')\left(1_{\bar{B}(c,r')}(x_{\xi'})-1_{\bar{B}(c,r')}(x_{\xi})\right)
\,\Pi(dX)\geq 0.
\end{equation}
Letting $r'\goto r+$ we have
\begin{equation}\label{equ6}
\begin{split}
&\int\sum_{\xi'\neq \xi}a(\xi,\xi')\left[1_{\bar{B}(c,r)}(x_\xi)1_{\bar{B}(c,r)}(x_{\xi'})-1_{\bar{B}(c,r)}(x_{\xi})\right] \Pi(dX)\\
&=\int 1_{\bar{B}(c,r)}(x_\xi)\sum_{\xi'\neq
\xi}a(\xi,\xi')\left(1_{\bar{B}(c,r)}(x_{\xi'})-1_{\bar{B}(c,r)}(x_{\xi})\right)
\Pi(dX)\geq 0.
\end{split}
\end{equation}
Since
\[1_{\bar{B}(c,r)}(x_\xi)1_{\bar{B}(c,r)}(x_{\xi'})-1_{\bar{B}(c,r)}(x_{\xi})\leq 0,\]
it follows from (\ref{equ6}) that for any $\xi'\in I_\xi$
\[1_{\bar{B}(c,r)}(x_\xi)1_{\bar{B}(c,r)}(x_{\xi'})=1_{\bar{B}(c,r)}(x_{\xi}),\]
$\Pi$ almost everywhere. Because $c$ and $r$ are arbitrary and $E$
is separable, we have $x_\xi=x_{\xi'}$ $\Pi$ almost everywhere.

\section{Reversibility} \setcounter{equation}{0}

Let ${\mathcal L}$ denote the generator of the Fleming-Viot process
with mutation, selection, recombination, and no migration on each
colony.

\begin{definition}\label{irre}
A generator $A$ is said to be irreducible if for all $x$ in $E$ and
any non-negative, non-zero measurable function $g\in C(E)$, there
exists $t>0$ such that $(P_t g)(x)>0$, where $P_t$ is the semigroup
generated by $A$.
\end{definition}

\begin{theorem}\label{le7}
Assume that there is no migration, and the mutation generator $A$ is
irreducible. Let $\Pi$ be the reversible measure for ${\mathcal
L}_{s,r,0}$. Then for each $\xi$ in $I$, \be
\label{fullsupport}\Pi\{X\in M_1(E)^I
:\text{supp}(X_{\xi})=E\}=1.\ee

A probability measure $\Pi$ in $M_1(M_1(E)^I)$ is reversible with
respect to ${\cal L}_{s,r,0}$ if and only if there are $\theta
>0$, $\mu$ in $M_1(E)$, and $h$ in $C(E)$ such that, for any $g$ in
$C(E)$, the mutation generator $A$ and recombination kernel
$\eta(x,y;dz)$ satisfy
\beq
&& Ag(x) + r\left[\int g(z)\eta(x,x;dz)-g(x)\right] =\frac{\theta}{2} [\langle \mu, g\rangle-g(x)],\\
&& \eta(x,y;dz)=\frac{1}{2}\left(\eta(x,x;dz)+\eta(y,y;dz)\right)
+(h(x)-h(y))(\delta_{x}(dz)-\delta_y(dz)). \eeq
\end{theorem}
\proof When there is no migration, the interacting system becomes a
system of independent Fleming-Viot processes. The theorem is then a
direct result of Proposition 3.1 and Theorem 2.2 in \cite{handa02}.
\eproof

\begin{theorem}\label{tt1}
Assume that $\rho >0$ and that $E$ is not a one point space. If the
mutation operator $A$ is irreducible, then there is no reversible
measure with respect to ${\cal L}_{s,r,\rho}$.
\end{theorem}
\proof Let $\Pi$ be reversible for ${\cal L}_{s,r,\rho}$. For any
$\xi$ in $\hat{I}$, (\ref{dirac}) shows that $X_\xi$ is $\Pi$-almost
surely a Dirac measure. On the other hand, the projection of $\Pi$
to each colony $\xi$ in $\hat{I}$ is a reversible measure of the
Fleming-Viot process on the colony. Applying Proposition 3.1 in
\cite{handa02} again it follows that $X_\xi$ has full support
$\Pi$-almost surely. This implies that $E$ is a one point space. A
contradiction. \eproof

We now consider the case of zero mutation. For any $\xi, \xi'\in I$,
write $\xi'\goto\xi$ if either $a(\xi',\xi)>0$ or there exists a
finite sequence $\xi_i, i=1,\ldots,n$ such that $a(\xi',\xi_1)>0,
a(\xi_1,\xi_2)>0, \ldots, a(\xi_n,\xi)>0$. Recall that  $\Delta$
denotes the collection of Dirac measures on $E$. Set
\[\Delta_a:=\{X\in M_1(E)^I: X_\xi= X_{\xi'}\in \Delta, \forall \xi, \xi'\in I \text{\, with \, }\xi'\goto \xi\}.\]

\begin{theorem}\label{rev}
Suppose that $\rho >0$ and for any $\xi$ in $I$, there is $\xi'$
such that $\xi'\ra \xi$. If there is no mutation or recombination,
then $\Pi$ is reversible if and only if its support is in
$\Delta_a$.
\end{theorem}

\proof The necessity follows from Lemma \ref{le5} and Lemma
\ref{le6}. If the mutation and recombination are zero, then for any
$X \in\Delta_a$, and any $F\in {\cal A}$ we have ${\mathcal
L}_{s,r,\rho}F(X) = 0$, which gives the sufficiency. \eproof

\section{Examples} \setcounter{equation}{0}
In this section, we discuss the reversibility of several well-known
examples.

{\bf Example 1.} (Two Type Stepping-Stone Model). Let $I= \IZ^d$ be
the $d$ dimensional lattice, and $E=\{0,1\}$. Let $x_i$ denote the
the proportion of type $0$ individuals on colony $i$ in $\IZ^d$. The
generator on colony $i$ is given by
\[L= \frac{1}{2}a_i(x)\frac{\partial^2 }{\partial x_i^2} +b_i(x)\frac{\partial }{\partial x_i},\]
where \beq
&& x=(x_i: i \in \IZ^d),\ a_i(x)= x_i(1-x_i),\\
&& b_i(x)=\sum_{j\in \IZ^d}\alpha(i,j)(x_j-x_i) + v-(u+v)x_i +s x_i(1-x_i), \\
&& \alpha(i,j)\geq 0, u,v \geq 0. \eeq This is the model studied in
\cite{shiga80} and \cite{shigauchi86}. It follows from
Theorem~\ref{tt1} and Theorem~\ref{rev} that the model has a
reversible measure if and only if $d\leq 2$ and $u=v=0$.

{\bf Example 2.} This model, studied in \cite{dgv95}, has zero
mutation and recombination. Let $I$ be either $\IZ^d$ or the
hierarchical group $\Omega_N$. In addition to assumption \rf{mig1},
the migration rate satisfies $a(\xi,\xi')=a(0,\xi'-\xi)$ and
$\sum_{n=0}^\infty\left(a^n(0,\xi)+a^n(\xi,0)\right)>0$. Set
$\hat{a}(\xi,\xi')=\frac{1}{2}[a(\xi,\xi')+a(\xi',\xi)]$.
 Theorem~\ref{tt1} and Theorem~\ref{rev} imply that the model has a reversible measure if and only if
the symmetrized kernel $\hat{a}$ is recurrent.

\section{Appendix}

 \noindent {\bf Definition.} Let $S$ be a metric space. A sequence
$\{h_n\}\subset B(S)$ is said to converge boundedly and pointwise to
$f\in B(S)$ if $h_n(x)\goto h(x)$ for all $x\in S$ and $\sup_n
\|h_n\|_\infty<\infty $. We write
\[{\text{bp}-\lim}_{n\goto\infty}h_n=h.\]

\noindent Part 1. The space $\cal H$.

Define $\cal H$ to be the space of functions $F:{\p}^I\to\IR$ so
that the partial derivative  $\de F(X)/\de X_\xi (u)$ exists for
every $X$, $\xi$, and $u$, and (\ref{t2b}) holds with $\Phi$
replaced by $F$.

Our first observation is that for any positive integer $m$, any
${\bf f}\in B(E^m)$ and any $(\xi_1,\ldots,\xi_m)\subset I^m$, the
function $F_{\bf f}: M_1(E)^I \mapsto \IR $ defined by $F_{\bf
f}(X):=\lan\otimes_{i=1}^m X_{\xi_i},\bf f \ran $ belongs to $\cal
H$. First consider the case of ${\bf f}= 1_{G_1\times\cdots\times
G_m}$ for open sets $G_i\subset E, i=1,\cdots,m$. Since we can
approximate the indicator function $1_{G_i}$ boundedly and pointwise
by functions in $D(A)$ which is dense in $C(E)$, it follows that one
can find a sequence of functions ${\bf f}_n$ in ${\cal A}$ such that
${\text{bp}-\lim}_{n\goto\infty}{\bf f}_n={\bf f}$. Since the
bp-convergence of ${\bf f}_n$ to ${\bf f}$ implies the
bp-convergence of the corresponding derivatives, we have that
$\lan\otimes_{i=1}^m X_{\xi_i},{\bf f} \ran\in\cal H$. Then the
observation follows from Theorem 4.3 in the Appendixes of
\cite{EtKu86}.

Using the above-mentioned observation and polynomial approximation
we can further show that for any $m_i$, any
$(\xi_{i1},\ldots,\xi_{im_i})\in I^{m_i}$, any ${\bf f}_i\in
B(\IR^{m_i}), i=1,\ldots,n$, and any $\phi\in C^1(R^n)$, the
function $F :M_1(E)^I \mapsto \IR$ defined by
\[F(X):=\phi(\lan\otimes_{j=1}^{m_1} X_{\xi_{1j}},{\bf
f}_1\ran,\ldots,\lan\otimes_{j=1}^{m_n} X_{\xi_{nj}},{\bf
f}_n\ran)\] also belongs to $\cal H$.

Moreover, take $g=\otimes_{i=1}^m g_i$ with $g_i\in {\cal D}(A)$
bounded below by $c>0$, and $k\in B(E^m)$ and set
$F(X):=\Phi^k(X)/\Phi^g(X)$. By polynomial approximation again we
can show that $F\in \cal H$.

\vskip 5pt \noindent Part 2. Approximating $\tilde{\Phi}_t$. \qquad

Let ${\bf f}\in {\cal D}(A)^I_0$ such that outside the finite subset
$I_0$ of $I$ $f_{\xi}\equiv 0$, and $X_s=S_{-s{\bf f}}X$. Then
\begin{equation*}
\begin{split}
\l b(X_s),{\bf f}\r & =\sum_{\xi\in I_0} \langle
X_{\xi}^{sf_\xi},Af_{\xi}\rangle
+ \rho\sum_{\xi\in I_0}\sum_{\xi' \in I}a(\xi,\xi') \langle X_{\xi'}^{sf_{\xi'}}-X_{\xi}^{sf_\xi}, f_{\xi}\rangle \\
&\quad + s\sum_{\xi \in I_0}\left(\int_E\int_E
V(u,v)f_\xi(u)X_\xi^{sf_\xi}(dv)X_\xi^{sf_\xi}(du)-\langle f_{\xi}, X_\xi^{sf_\xi}\rangle \langle V, {X_\xi^{sf_\xi}}^{\otimes 2}\rangle\right)\\
&\quad+r\sum_{\xi \in I_0}\left(\left\langle \int_E
f_\xi(u)\eta(\cdot,\cdot;du), {X_\xi^{sf_\xi}}^{\otimes
2}\right\rangle-\langle f_{\xi}, X_\xi^{sf_\xi}\rangle\right).
\end{split}
\end{equation*}
Since $\sum_{\xi' \in I}a(\xi,\xi')<\infty$ and
\[\sum_{\xi\in I_0}\sum_{\xi' \in I}a(\xi,\xi') \langle X_{\xi'}^{sf_{\xi'}}, f_{\xi}\rangle
=\sum_{\xi\in I_0}\sum_{\xi' \in I_0}a(\xi,\xi') \langle
X_{\xi'}^{sf_{\xi'}}, f_{\xi}\rangle+ \sum_{\xi\in I_0}\sum_{\xi'
\not\in I_0}a(\xi,\xi') \langle X_{\xi'}, f_{\xi}\rangle,\] by Part
1 we have $\l b(X_s),{\bf f}\r \in {\cal H}$.

Define $\Phi_t(X):=\Phi(X_t)$,
\[\Lambda_n({\bf f},X):=\frac{2t}{n}\sum_{i=1}^n \l b(X_{it/n}), {\bf f}\r\]
and
\[\tilde{\Phi}_t^{(n)}(X):=\Phi(X_t)e^{-\Lambda_n({\bf f},X)}.\]
Since both $\Phi(X_t)\in {\cal H}$  and $e^{-\Lambda_n({\bf
f},X)}\in {\cal H}$ by Part 1, then $\tilde{\Phi}_t^{(n)}\in {\cal
H} $ and (\ref{t2b}) holds with $\Phi$ replaced by
$\tilde{\Phi}_t^{(n)}$.

Clearly,
\[{\text{bp}-\lim}_{n\goto\infty}\tilde{\Phi}_t^{(n)}=\tilde{\Phi}_t.\]
Similar to (\ref{t2a}), we have
\begin{equation*}
\begin{split}
{\de \tilde\Phi_t^{(n)}(X)\over\de X_\xi}(u) &= \frac{\de
\Phi(X_t)}{\de X_\xi}(u) e^{-\Lambda_n({\bf f},X)} -
2\tilde\Phi_t^{(n)}(X) \frac{t}{n}\sum_{i=1}^n {\de \l
b(X_{it/n}),{\bf f}\r \over\de X_\xi}(u).
\end{split}
\end{equation*}
Therefore, \[{\text{bp}-\lim}_{n\goto\infty}\frac{\de
\tilde\Phi_t^{(n)}}{\de X_\xi}=\frac{\de\tilde{\Phi}_t}{\de X_\xi},
\,\,\, \forall\, \xi\in I,\] and (\ref{t2b}) holds for
$\tilde{\Phi}_t $.

      \end{document}